\documentclass{article}

\usepackage[english]{babel}

\usepackage{amssymb}
\usepackage{amsmath}
\usepackage{comment}
\usepackage{secdot}
\newcounter{nRemark}

\renewcommand{\thenRemark}{\arabic{nRemark}}
\newenvironment{remark}
{\refstepcounter{nRemark}%
\addvspace{\smallskipamount}\noindent
\textbf{Remark~\thenRemark.~}} {\par\addvspace{\smallskipamount}}

\newenvironment{remark*}
{\addvspace{\smallskipamount}\noindent
\textbf{Remark}.~} {\par\addvspace{\smallskipamount}}
\newcounter{nTheorem}

\renewcommand{\thenTheorem}{\arabic{nTheorem}}
\newenvironment{Theorem}
{\refstepcounter{nTheorem}%
\addvspace{\smallskipamount}\noindent\leftskip1.5\parindent\parindent-\parindent
\textbf{Theorem~\thenTheorem.~}} {\par\addvspace{\smallskipamount}}
\newenvironment{proof}
{\addvspace{\smallskipamount}\noindent
\textsf{\textit{Proof.}}} {\par\addvspace{\smallskipamount}}

\begin{document}

\title{Affine transformations of the plane and their geometrical properties}
\author{Irina Busjatskaja \and Yury Kochetkov}
\date{}
\maketitle
\begin{abstract}
In this work, oriented for students with knowledge of basics of
linear algebra, we demonstrate, how the use of polar decomposition
allows one to understand metric properties of non-degenerate
linear operators in \(\mathbf{R}^2\).
\end{abstract}

\section*{Introduction}

Let \(\mathbf{R}^2\) be the standard two-dimensional Euclidean space and \(\varphi\) be a
reversible, non-isometric linear operator on this space. Let
\(A=\begin{pmatrix}
a & b \\
c & d \end{pmatrix}\) be the matrix of \(\varphi\) in the standard basis
\(\mathbf{i}=\begin{pmatrix}
1 \\
0 \end{pmatrix}\),
\(\mathbf{j}=\begin{pmatrix}
0 \\
1 \end{pmatrix}\) of \(\mathbf{R}^2\). We will assume that operator \(\varphi\) is reversible, so
\(\det(A)=ad-bc\neq0\).

Let \(P(\lambda)=\begin{vmatrix}
a-\lambda & b \\
c & d-\lambda \end{vmatrix}\) be the characteristic
polynomial of operator \(\varphi\). The roots \(\lambda_1,\lambda_2\)
of polynomial \(P\) could be either real numbers (real spectrum), or complex conjugate (complex
spectrum). In the case of real spectrum we will assume that
\(\lambda_2>\lambda_1>0\). In the real case there are two one-dimensional invariant subspaces
\(L_1\) and \(L_2\) in \(\mathbf{R}^2\), which consist of the eigenvectors with eigenvalues \(\lambda_1\)
and \(\lambda_2\), respectively. In the complex case there are no non-trivial
invariant subspaces.

Matrix \(A\) of operator \(\varphi\) satisfy the condition
\begin{equation}
 \label{eq1:art2}
(a+d)^{2}-4(ad-bc)<0
\end{equation}
or
\begin{equation}
 \label{eq2:art2}
(a-d)^{2}+4bc<0.
\end{equation}
in the complex case, and the condition
\begin{equation}
 \label{eq3:art2}
(a+d)^{2}-4(ad-bc)>0
\end{equation}
in the real case.

\begin{remark*}
As similar matrices have the same trace and determinant, then conditions
\eqref{eq1:art2}, \eqref{eq2:art2} and \eqref{eq3:art2} will be true for
matrix of operator in any base.

In this paper we will determine the range of angles between the vectors
\(\mathbf{x}\) and \(\varphi(\mathbf{x})\), using the polar decomposition of given operator \(\varphi\).
In \ref{s1:art2} we
will recall the notion of polar decomposition and its properties. In \ref{s2:art2}
we will consider an operator with a real spectrum and will estimate the
angle \(\gamma\) between the vectors \(\mathbf{x}\) and \(\varphi(\mathbf{x})\).
\end{remark*}

\begin{Theorem}
Let \(\varphi\) be the operator with the real spectrum, \(\lambda_1,\lambda_2\) are positive
eigenvalues of the operator \(\varphi\) and \(\beta\) is the acute angle between the
eigenvectors, then
\(0\leqslant\gamma\leqslant\arccos\dfrac{2\sqrt{\lambda_1\lambda_2}-(\lambda_1+\lambda_2)\cos\beta}{\lambda_1+\lambda_2-2\sqrt{\lambda_1\lambda_2}\cos\beta}\).
\end{Theorem}

In \ref{s3:art2} we consider an operator \(\varphi\) with the complex spectrum and,
using its polar decomposition, we obtain the following result:

\begin{Theorem}
Let \(\varphi\) be an operator with the complex spectrum, then
\(\alpha-\arccos\dfrac{2\sqrt[4]{\lambda\mu}}{\sqrt{\lambda}+\sqrt{\mu}}\leqslant\gamma\leqslant
\alpha+\arccos\dfrac{2\sqrt[4]{\lambda\mu}}{\sqrt{\lambda}+\sqrt{\mu}}\),
where \(\alpha\) is the rotation angle of the isometric operator and \(\sqrt{\lambda}\) and
\(\sqrt{\mu}\) are the eigenvalues of the positive self-conjugate operator in the
polar decomposition of the operator \(\varphi\).
\end{Theorem}

In \ref{s4:art2} we compute \(\bar{\gamma}\) the mean value of the angle \(\gamma\) for the
operators with the complex spectrum.

\begin{Theorem}
For the operators with the complex spectrum
\(\dfrac{\pi}{2}-\dfrac{2}{\pi}\leqslant\bar{\gamma}\leqslant\dfrac{\pi}{2}+\dfrac{2}{\pi}\).
\end{Theorem}

In \ref{s5:art2} we compute the norm of an operator \(\varphi\) and in
\ref{s6:art2} we study trajectories of the points under the action of an
operator with the complex spectrum such that \(\det(A)=1\). We prove that each
trajectory belongs to some ellipse.

Other applications of polar decomposition can be seen in \cite{4-art2}.

\section{The polar decomposition of the reversible linear operator on the plane}\label{s1:art2}
Let us note at first that if \(\varphi\) is an operator in \(\mathbf{R}^2\) with complex spectrum,
then \(\det(A)>0\), where
\(A=\begin{pmatrix}
a & b \\
c & d \end{pmatrix}\) is the matrix of \(\varphi\) in the
standard base. Indeed, in this case the discriminant \((a+d)^2-4(ad-bc)\) of the characteristic
polynomial of matrix \(A\) is negative.

In the real case we will assume that eigenvalues \(\lambda_1\) and \(\lambda_2\) of \(\varphi\) are
positive. Here \(\det(A)\) is also positive.

Now let us consider a self-conjugate operator
\(\varphi\circ\varphi^*\) (\(\varphi\varphi^*(x)=\varphi^*\bigl(\varphi(x)\bigr)\))
and its (symmetric) matrix \(A^{\mathrm{t}}A\). It is well known that eigenvalues \(\lambda\) and
\(\mu\) of operator \(\varphi\circ\varphi^*\) are positive and the corresponding eigenvectors \(\mathbf{e}_1\)
and \(\mathbf{e}_2\)~--- orthonormal. There exists operator \(\varphi^{\prime}\) such that
\(\varphi\circ\varphi^*=(\varphi^{\prime})^2\).
\(\varphi^{\prime}\) is a positive operator with the matrix \(B\) in the standard base and the
matrix
\begin{equation}
 \label{eq4:art2}
B^{\prime}=
\begin{pmatrix}
\sqrt{\lambda } & 0 \\
0 & \sqrt{\mu }
\end{pmatrix}
\end{equation}
in base \(\mathbf{e}_1\) and \(\mathbf{e}_2\). It must be noted that \(B\) is a symmetric matrix,
because \(\varphi^{\prime}\) is a self-adjoint operator, and \(B^2=A^{\mathrm{t}}A\).

The product \(A=(AB^{-1})B=OB\) is called a polar decomposition of matrix
\(A\). Here the matrix \(O=AB^{-1}\) is orthogonal. Indeed,
\[
O^{\mathrm{t}}O=(AB^{-1})^{\mathrm{t}} AB^{-1} = B^{-1} A^{\mathrm{t}}A B^{-1} = B^{-1} B^2 B^{-1} = E.
\]
Orthogonal matrix \(O\) defines an operator \(\psi\), which is either a rotation, if
\(\det(O)>0\), or reflection, if \(\det(O)<0\). If \(\det(A)>0\) and \(\varphi^{\prime}\) is a positive
operator, then \(O\) defines a rotation on angle \(\alpha\), thus
\[
O=
\begin{pmatrix}
\cos\alpha & -\sin\alpha \\
\sin\alpha & \cos\alpha
\end{pmatrix}.
\]
If \(\varphi\) is an operator with complex spectrum, then, as we know, \(\det(A)>0\).
So, in this case also the orthogonal operator \(\psi\) in polar decomposition is a rotation.

Let \(A^{\prime}\) be the matrix of operator \(\varphi\) in the base \(\mathbf{e}_1\) and \(\mathbf{e}_2\).
As \(A^{\prime}= O B^{\prime}\), then
\[
A^{\prime}=
\begin{pmatrix}
\cos\alpha & -\sin\alpha \\
\sin\alpha & \cos\alpha
\end{pmatrix}
\begin{pmatrix}
\sqrt{\lambda } & 0 \\
0 & \sqrt{\mu }
\end{pmatrix}.
\]
Condition \eqref{eq2:art2} may be rewritten as
\[
(\sqrt{\lambda }\cos\alpha - \sqrt{\mu }\cos\alpha)^{2} - 4\sqrt{\lambda }\,\sqrt{\mu }\sin^{2}\alpha<0.
\]
Thus
\begin{equation}
 \label{eq5:art2}
\cos\alpha<\frac{2\sqrt[4]{\lambda\mu}}{\sqrt{\lambda}+\sqrt{\mu}}
\end{equation}
in the case of complex spectrum, and
\begin{equation}
 \label{eq6:art2}
\cos\alpha>\frac{2\sqrt[4]{\lambda\mu}}{\sqrt{\lambda}+\sqrt{\mu}}
\end{equation}
in the case of the real one.

\section{Angles of rotations of the operator with the real spectrum} \label{s2:art2}

Let \(\lambda_1\) and \(\lambda_2\), \(0<\lambda_1<\lambda_2\), be eigenvalues
of \(\varphi\) and \(\mathbf{u}_1\) and \(\mathbf{u}_2\),
\(\lvert\mathbf{u}_1\rvert=\lvert\mathbf{u}_2\rvert=1\), be the corresponding
eigenvectors.
Let \(\beta\) be the acute angle between \(\mathbf{u}_1\) and \(\mathbf{u}_2\).
The plane is divided in four parts (cones) by two lines corresponding to
the subspaces \(L_1 = \langle\mathbf{u}_1\rangle\) and \(L_2 = \langle\mathbf{u}_2\rangle\).
Each open cone is invariant under the action of \(\varphi\) and
under this action each vector \(\mathbf{x}\) moves in the direction to the line
\(L_2\). Let us compute the maximal value \(\gamma_{\text{max}}\) of the angle \(\gamma\)
between the vectors \(\mathbf{x}\) and \(\varphi(\mathbf{x})\)  and estimate the angle \(\gamma\).

\setcounter{nTheorem}{0}
\begin{Theorem}
Let \(\varphi\) be the operator with the real spectrum , \(\lambda_1, \lambda_2\) are positive
eigenvalues of the operator \(\varphi\) and \(\beta\) is the acute angle between the
eigenvectors, then \(0 \leqslant \gamma \leqslant
\arccos\dfrac{2\sqrt{\lambda_1\lambda_2}-(\lambda_1+\lambda_2)\cos\beta}{\lambda_1+\lambda_2-2\sqrt{\lambda_1\lambda_2}\cos\beta}\).
\end{Theorem}

\begin{proof}
We shall work in the base \(\{\mathbf{u}_1,\mathbf{u}_2\}\). Let \(x_1\) and \(x_2\) be the
coordinates of vector \(\mathbf{x}\) in our base, then
\(\varphi(\mathbf{x})= \lambda_1 x_1 \mathbf{u}_1+\lambda_2 x_2 \mathbf{u}_2\). As
\(\cos\gamma=\dfrac{(\mathbf{x},\varphi(\mathbf{x}))}{\lvert \mathbf{x}\rvert\,\lvert\varphi(\mathbf{x})\rvert}\)
and Gram matrix of the vectors \(\mathbf{u}_1,\mathbf{u}_2\) is
\(G=
\begin{pmatrix}
1 & \cos\beta \\
\cos\beta & 1
\end{pmatrix}\) then
\(\cos\gamma =\dfrac{(x_1 x_2)\,
G
\begin{pmatrix}
\lambda_1 & x_1 \\
\lambda_2 & x_2
\end{pmatrix}}{\sqrt{(x_1\,x_2)\,
G
\begin{pmatrix}
x_1 \\
x_2
\end{pmatrix}
(\lambda_1x_1\,\lambda_2x_2)
G
\begin{pmatrix}
\lambda_1 & x_1 \\
\lambda_2 & x_2
\end{pmatrix}}}\).
\begin{equation}
 \label{eq7:art2}
\cos\gamma = \dfrac{\lambda_1(x_1)^2 + (\lambda_1 + \lambda_2)x_1x_2\cos\beta + \lambda_2(x_2)^2}
{\sqrt{x_1 + 2x_1x_2\cos\beta + (x_2)^2}\,
\sqrt{(\lambda_1)^2(x_2)^2 + 2\lambda_1 \lambda_2 x_1x_2\cos\beta + (\lambda_2)^2(x_2)^2}}
\end{equation}
If we denote \(\dfrac{x_1}{x_2}\) by \(t\) and \(\dfrac{\lambda_2}{\lambda_1}\) by \(\delta\),
then the formula \eqref{eq7:art2} may be rewritten as
\begin{equation}
 \label{eq8:art2}
f(t) = \cos\gamma = \frac{(t)^2 + (1+\delta)\cos\beta t + \delta}
{\sqrt{(t)^2 + 2t\cos\beta + 1}\,\sqrt{(t)^2 + 2\delta t\cos\beta + (\delta)^2}}.
\end{equation}
To estimate the values of the function \(f(t)\) we calculate its derivative
\begin{equation}
 \label{eq9:art2}
f^{\prime}(t) = \frac{(\delta-1)^2 (1 - (\cos\beta)^2) ((t)^3 - \delta t)}
{\sqrt{((t)^2 + 2t\cos\beta + 1)^3}\,\sqrt{((t)^2 + 2\delta t\cos\beta + (\delta)^2)^3}}.
\end{equation}
From \eqref{eq9:art2} follows that the equation \(f^{\prime}(t)=0\) has three roots: \(t=0\), \(t =\sqrt{\delta}\)
and \(t=-\sqrt{\delta}\), hence \(f_{\text{min}}(t) = f(\sqrt{\delta}) =
\dfrac{2\sqrt{\delta} + (1+\delta) \cos\beta}{1+2\sqrt{\delta}\cos\beta + \delta}\) for
\(t\in (0,+\infty)\), and
\(f_{\text{min}}(t) = f(-\sqrt{\delta}) =
\dfrac{2\sqrt{\delta} - (1+\delta) \cos\beta}{1 - 2\sqrt{\delta}\cos\beta + \delta}\)
for \(t\in (-\infty,0)\). Notice that if \(\cos\beta\neq0\) then \(f(-\sqrt{\delta}) < 0 < f(\sqrt{\delta})\)
and hence
\begin{equation}
 \label{eq10:art2}
\gamma_{\text{max}} = \arccos\frac{2\sqrt{\lambda_1\lambda_2}-(\lambda_1+\lambda_2)\cos\beta}
{\lambda_1+\lambda_2-2\sqrt{\lambda_1\lambda_2}\cos\beta}.
\end{equation}

Thus \(0 \leqslant \gamma \leqslant \arccos\dfrac{2\sqrt{\lambda_1\lambda_2}-(\lambda_1+\lambda_2)\cos\beta}
{\lambda_1+\lambda_2-2\sqrt{\lambda_1\lambda_2}\cos\beta}\).
\end{proof}

\begin{remark}
If the base \(\mathbf{u}_1,\mathbf{u}_2\) is orthonormal, thus \(\cos\beta =0\), then from
\eqref{eq10:art2} follows that
\begin{equation}
 \label{eq11:art2}
\gamma_{\text{max}} = \arccos\frac{2\sqrt{\lambda_1\lambda_2}}{\lambda_1+\lambda_2}.
\end{equation}
\end{remark}

\begin{remark}
If \(\lambda_1<0\), \(\lambda_2<0\), then the operator \(\varphi\) is the composition of the
operator with the real positive spectrum \(\lvert\lambda_1\rvert,\lvert\lambda_2\rvert\) and the operator
of the central symmetry and angle \(\gamma\) also can be estimated
\(\pi - \arccos\dfrac{2\sqrt{\lambda_1\lambda_2}-(\lambda_1+\lambda_2)\cos\beta}
{\lambda_1+\lambda_2-2\sqrt{\lambda_1\lambda_2}\cos\beta} \leqslant \gamma \leqslant \pi\).
\end{remark}

\begin{remark}
If \(\lambda_1\) and \(\lambda_2\) are of different signs, then vectors
\(\mathbf{x}\) and \(\varphi(\mathbf{x})\) are in the adjacent cones and
\(0\leqslant\gamma\leqslant\pi\).
\end{remark}

\section{Angles of rotations of the operator with the complex spectrum} \label{s3:art2}

Let
\begin{equation}
\label{eq12:art2}
A^{\prime}=
\begin{pmatrix}
\cos\alpha & -\sin\alpha \\
\sin\alpha & \cos\alpha
\end{pmatrix}
\begin{pmatrix}
\sqrt{\lambda } & 0 \\
0 & \sqrt{\mu }
\end{pmatrix}\quad\alpha\neq k\pi
\end{equation}
be its matrix in the base \(\mathbf{e}_1,\mathbf{e}_2\) (see \ref{s1:art2}).

Here \(\alpha\) is the rotation angle of the orthogonal operator \(\psi\) and
\(\lambda\) and \(\mu\) are eigenvalues of positive operator
\(\varphi\circ\varphi^* = (\varphi^{\prime})^2\). As \(\varphi\) has no eigenvectors
then it rotates each vector \(\mathbf{x}\neq0\) in one direction on some
angle \(\gamma=\gamma(\mathbf{x})\neq0\). If \(\gamma^{\prime}(\mathbf{x})\) is the angle
between vectors \(\mathbf{x}\) and \(\varphi^{\prime}(\mathbf{x})\),
then \(\gamma(\mathbf{x})=\alpha + \gamma^{\prime}(\mathbf{x})\). It must be noted
that angles \(\gamma^{\prime}\) may be positive and
negative i.\,e. rotation from \(\mathbf{x}\) to \(\varphi^{\prime}(\mathbf{x})\)
may be counterclockwise or clockwise.

\begin{Theorem}
Let \(\varphi\) be the operator with the complex spectrum, then
\[
\alpha-\arccos\frac{2\sqrt[4]{\lambda\mu}}{\sqrt{\lambda}+\sqrt{\mu}}\leqslant\gamma\leqslant
\alpha+\arccos\frac{2\sqrt[4]{\lambda\mu}}{\sqrt{\lambda}+\sqrt{\mu}}.
\]
\end{Theorem}

\begin{proof}
As \(\varphi^{\prime}\) is an operator with real positive spectrum and the base
\(\mathbf{e}_1,\mathbf{e}_2\) is orthonormal, then we can use the formula \eqref{eq11:art2}.
\begin{equation}
 \label{eq13:art2}
\gamma_{\text{max}}^{\prime} = \arccos\frac{2\sqrt[4]{\lambda\mu}}{\sqrt{\lambda}+\sqrt{\mu}}.
\end{equation}
Now from \eqref{eq5:art2} we see that \(\lvert\gamma_{\text{max}}^{\prime}\rvert < \lvert\alpha\rvert\),
thus maximal rotation angle of \(\varphi\) is the sum
\(\gamma_{\text{max}} = \alpha + \arccos\dfrac{2\sqrt[4]{\lambda\mu}}{\sqrt{\lambda}+\sqrt{\mu}}\) and minimal
rotation angle of \(\varphi\) is the difference
\(\gamma_{\text{min}} = \alpha - \arccos\dfrac{2\sqrt[4]{\lambda\mu}}{\sqrt{\lambda}+\sqrt{\mu}}\). So
\(\alpha-\arccos\dfrac{2\sqrt[4]{\lambda\mu}}{\sqrt{\lambda}+\sqrt{\mu}}\leqslant\gamma\leqslant
\alpha+\arccos\dfrac{2\sqrt[4]{\lambda\mu}}{\sqrt{\lambda}+\sqrt{\mu}}\).
\end{proof}

\begin{remark*}
For an operator \(\varphi\) with the real positive spectrum, its polar
decomposition also can be used to calculate the angles \(\gamma_{\text{min}}\)
and \(\gamma_{\text{max}}\), but the condition \eqref{eq6:art2} shows us
that \(\lvert\gamma_{\text{max}}\rvert > \lvert\alpha\rvert\), hence \(\gamma_{\text{min}}\)
and \(\gamma_{\text{max}}\) have different signs. It means that \(\varphi\) rotates vectors in
different directions, while the operator \(\varphi\) with the complex spectrum
rotates all vectors in one direction~--- in the direction of the angle \(\alpha\).
\end{remark*}

\section{The mean value of the angle of rotation (complex case)}\label{s4:art2}

Let \(\varphi\) be an operator with the complex spectrum and
\(A=
\begin{pmatrix}
a & b \\
c & d
\end{pmatrix}\)
be its matrix in the standard base. Then, as we know \eqref{eq1:art2},
\(a^2+2ad+4bc+d^2<0\). This condition defines a domain \(G\) in the
four-dimensional real space with coordinates \(a, b ,c, d\).

In \ref{s3:art2} it was proved that \(\alpha - \gamma_{\text{max}}^{\prime} < \gamma < \alpha + \gamma_{\text{max}}^{\prime}\)
where \(\gamma\) is the angle between the vectors \(\mathbf{x}\) and \(\varphi(\mathbf{x})\),
\(\alpha\) is the rotation angle of orthogonal operator \(\psi\) and
\(\gamma_{\text{max}}^{\prime} = \arccos\dfrac{2\sqrt[4]{\lambda\mu}}{\sqrt{\lambda}+\sqrt{\mu}}\)
is the maximal rotation angle of positive operator \(\varphi^{\prime}\) (\(\varphi\circ\varphi^*=(\varphi^{\prime})^2\)).
As \(\lambda\) and \(\mu\) are eigenvalues of matrix \(A^{\mathrm{t}}A\), then \(\lambda + \mu = a^2+b^2+c^2+d^2\),
\(\lambda\mu = (ad-bc)^2\). Set \(a=\dfrac{1}{\sqrt{2}}(x + t)\), \(b=\dfrac{1}{\sqrt{2}}(y + z)\),
\(c = \dfrac{1}{\sqrt{2}}(y - z)\), \(d = \dfrac{1}{\sqrt{2}}(t - x)\), then
\(\sqrt{\lambda}\,\sqrt{\mu} = \dfrac12(z^2 + t^2 - x^2 - y^2)\),
\(\sqrt{\lambda} + \sqrt{\mu} = \sqrt{2(t)^2+2(z)^2}\) and the inequality \(x^2 + y^2 - z^2 < 0\)
describes the domain \(G\) in new coordinates.

Also we have, that \(\gamma_{\text{max}}^{\prime} = \gamma_{\text{max}}^{\prime}(x,y,z,t) =
\arccos\sqrt{1 - \dfrac{(x)^2+(y)^2}{(t)^2+(z)^2}} = \arcsin\sqrt{\dfrac{(x)^2+(y)^2}{(t)^2+(z)^2}}\).
Let \(\bar{\gamma}\) be the mean value of the angle of rotation of the operators
with the complex spectrum.

\begin{Theorem}
For the operators with the complex spectrum
\(\dfrac{\pi}{2}-\dfrac{2}{\pi}\leqslant\bar{\gamma}\leqslant\dfrac{\pi}{2}+\dfrac{2}{\pi}\).
\end{Theorem}

\begin{proof}
As \(\alpha - \gamma_{\text{max}}^{\prime} < \gamma < \alpha + \gamma_{\text{max}}^{\prime}\), then,
in order to estimate \(\bar{\gamma}\), we have to calculate \(\bar{\gamma}_{\text{max}}^{\prime}\)~---
the mean value of the maximal rotation angle of operators \(\varphi^{\prime}\).
\begin{equation}
 \label{eq14:art2}
\bar{\gamma}_{\text{max}}^{\prime} = \frac{\int_{G} \gamma_{\text{max}}^{\prime}(x,y,z,t)\,dv}{\int_G dv}.
\end{equation}

Let us consider the domain \(G^{\prime} = G \cap C\), where \(C = C^1\times C^2\),
\(C^1=\{(x,y); (x)^2 + (y)^2 < 1\}\) and, \(C^2=\{(z,t); (z)^2+(t)^2<1\}\).
As \(\gamma_{\text{max}}^{\prime}(x,y,z,t)\) is a homogeneous function,
we , instead of calculating \(\bar{\gamma}_{\text{max}}^{\prime}\)
by the formula \eqref{eq14:art2}, can use the following one:
\begin{equation}
 \label{eq15:art2}
\bar{\gamma}_{\text{max}}^{\prime} =
\frac{\int_{G^{\prime}} \gamma_{\text{max}}^{\prime}(x,y,z,t)\,dv}{\int_{G^{\prime}} dv}.
\end{equation}

In order to simplify the integration we apply the polar coordinates
\(x =\rho\cos\eta\), \(y=\rho\sin\eta\), \(t=r\cos\xi\), \(z=r\sin\xi\).
\(G^{\prime}=\{(\eta,\xi,r,\rho);
0\leqslant\eta\leqslant2\pi,
0\leqslant\xi\leqslant2\pi,
0\leqslant r\leqslant1,
0\leqslant\rho\leqslant\lvert r\sin\xi\rvert\}\). Now we can compute both
integrals in formula \eqref{eq15:art2}.

\begin{align}
 \label{eq16:art2}
\int_{G^{\prime}}dv &= 4\int_0^{2\pi} d\eta \int_0^1 r\,dr \int_0^{\tfrac{\pi}{2}} d\xi
\int_0^{r\sin\xi} \rho\,d\rho = \frac{(\pi)^2}{4}\\
 \label{eq17:art2}
\int_{G^{\prime}} \gamma_{\text{max}}^{\prime}\,dv &= 4\int_0^{2\pi} d\eta \int_0^1 r\,dr \int_0^{\tfrac{\pi}{2}} d\xi
\int_0^{r\sin\xi} \rho\arcsin\frac{\rho}{r}\,d\rho = \frac{\pi}{8}.
\end{align}
So from \eqref{eq15:art2}--\eqref{eq17:art2} it follows that \(\bar{\gamma}_{\text{max}}^{\prime}=\dfrac{2}{\pi}\).

In order to get the sensible result for \(\bar{\gamma}\) we have to
limit the value of the angles \(\alpha\)~--- the angles of rotation of the
operators \(\psi\). Let \(\alpha\in[0,\pi]\) then \(\alpha\) can be uniquely
determined by \(\cos\alpha\). As similar matrices have the same trace,
then \(\cos\alpha (\sqrt{\lambda}+\sqrt{\mu}) = a + d\) or, after the
change of variables, \(\cos\alpha = \dfrac{t}{\sqrt{(t)^2+(z)^2}}\). As
\(\sin\alpha> 0\) for \(\alpha\in[0,\pi]\) then \(\sin\alpha =
\dfrac{z}{\sqrt{(t)^2+(z)^2}}\), \(\alpha = \arccos\dfrac{t}{\sqrt{(t)^2+(z)^2}}\)
and, in order to compute \(\bar{\alpha}\)~--- the mean value of the
angles $\alpha$, we have to consider only part of the domain
\(G\) corresponding to \(z>0\). As \(\alpha = \alpha(x,y,z,t)\) is the homogeneous
function, we can use the domain \(G^{\prime}\) with the condition \(z>0\).
\[
\bar{\alpha} = \frac{2}{(\pi)^2}\int_0^{2\pi} d\eta \int_0^1 r\,dr \int_0^{\pi}\xi\,d\xi
\int_0^{r\sin\xi} \rho\,d\rho = \frac{\pi}{2}.
\]
Thus we have the final result
\(\dfrac{\pi}{2}-\dfrac{2}{\pi}\leqslant\bar{\gamma}\leqslant\dfrac{\pi}{2}+\dfrac{2}{\pi}\).
\end{proof}

\section{The norm of the operator}\label{s5:art2}

Let now \(\varphi\) be an operator with the complex spectrum or an operator with
the positive real spectrum, then \(A^{\prime}=
\begin{pmatrix}
\cos\alpha & -\sin\alpha \\
\sin\alpha & \cos\alpha
\end{pmatrix}
\begin{pmatrix}
\sqrt{\lambda } & 0 \\
0 & \sqrt{\mu }
\end{pmatrix}\quad\alpha\neq k\pi\)
be its matrix in the base \(\mathbf{e}_1,\mathbf{e}_2\) (see \ref{s1:art2}).
The operator \(\varphi\) not only rotates the vector \(\mathbf{x}\) but also changes
its length. If \(\mathbf{x}=x_1 \mathbf{e}_1 + x_2 \mathbf{e}_2\) then
\(\varphi^{\prime}(\mathbf{x})= \sqrt{\lambda}\,x_1 \mathbf{e}_1 + \sqrt{\mu}\,x_2 \mathbf{e}_2\).
As an orthogonal operator \(\psi\)
does not change the length of the vectors, then
\(\lvert\varphi(\mathbf{x})\rvert = \lvert\varphi^{\prime}(\mathbf{x})\rvert\).
So \(\dfrac{(\lvert\varphi(\mathbf{x})\rvert)^2}{(\lvert \mathbf{x}\rvert)^2} = \dfrac{\lambda(x_1)^2 + \mu(x_2)^2}{(x_1)^2+(x_2)^2} =
\mu+\dfrac{(\lambda-\mu)(x_1)^2}{(x_1)^2+(x_2)^2}\). It is easy to estimate all
possible values of this ratio and to calculate the norm of an operator \(\varphi\).
\begin{equation}
 \label{eq18:art2}
\min(\lambda,\mu)\leqslant\frac{(\lvert\varphi(\mathbf{x})\rvert)^2}{(\lvert \mathbf{x}\rvert)^2}\leqslant\max(\lambda,\mu).
\end{equation}

If \(x_1=0\), then \(\dfrac{\lvert\varphi(\mathbf{x})\rvert}{\lvert \mathbf{x}\rvert}=\sqrt{\mu}\),
but if \(x_2=0\), then
\(\dfrac{\lvert\varphi(\mathbf{x})\rvert}{\lvert \mathbf{x}\rvert}=\sqrt{\lambda}\),
hence from \eqref{eq18:art2} it follows that the norm of an
operator \(\varphi\) is \(\|\varphi\|=\max(\sqrt{\lambda},\sqrt{\mu})\).

Notice that if \(\dfrac{(x_2)^2}{(x_1)^2} = \dfrac{\lambda-1}{1-\mu}\) then an operator
\(\varphi\) preserves the length of the vector \(\mathbf{x}\). We can find vector \(\mathbf{x}\)
satisfying this condition if and only if \(\min(\lambda,\mu)\leqslant1\leqslant\max(\lambda,\mu)\).

\section{Some curves related to the operator with the complex spectrum}\label{s6:art2}

Let \(\varphi\) be an operator with the complex spectrum, such that
\(\det(A)=1\). Thus its eigenvalues are complex conjugate numbers
\(\exp(i\theta)\), \(\exp(-i\theta)\) of modul \(1\).

Let \(\Phi\) be the complexification of \(\varphi\) and
\(\mathbf{z}=\begin{pmatrix}
z_1 \\
z_2
\end{pmatrix}
=\begin{pmatrix}
u_1 \\
u_2
\end{pmatrix}
+ i
\begin{pmatrix}
v_1 \\
v_2
\end{pmatrix} = \mathbf{u} + i\mathbf{v}\)
be the eigenvector of \(\Phi\) corresponding to the eigenvalue \(\exp(i\theta)\),
then \(\Phi(\mathbf{z}) = A\cdot
\begin{pmatrix}
z_1 \\
z_2
\end{pmatrix}
= A \cdot
\begin{pmatrix}
u_1 \\
u_2
\end{pmatrix}
+ i A \cdot
\begin{pmatrix}
v_1 \\
v_2
\end{pmatrix}
= \varphi(\mathbf{u}) + i\varphi(\mathbf{v})\), where \(\mathbf{u}\in\mathbf{R}^2\),
\(\mathbf{v}\in\mathbf{R}^2\).

As \(\Phi(\mathbf{z})=\exp(i\theta)\cdot\mathbf{z} = (\cos\theta + i \sin\theta)\cdot(\mathbf{u}+ i \mathbf{v}) =
(\cos\theta\cdot\mathbf{u} - \sin\theta\cdot\mathbf{v}) + i(\sin\theta\cdot\mathbf{u} + \cos\theta\cdot\mathbf{v})\),
then
\begin{align}
 \label{eq19:art2}
\varphi(\mathbf{u}) &= \cos\theta \cdot \mathbf{u} - \sin\theta \cdot \mathbf{v}, \\
 \label{eq20:art2}
\varphi(\mathbf{v}) &= \sin\theta \cdot \mathbf{u} + \cos\theta \cdot \mathbf{v}.
\end{align}

Now consider conjugate complex vector \(\bar{\mathbf{z}} = \mathbf{u} - i\mathbf{v}\).
\(\Phi(\bar{\mathbf{z}}) = A \cdot \bar{\mathbf{z}} = \overline{A \cdot \mathbf{z}} =
\overline{\exp(i\theta)} \cdot \bar{\mathbf{z}} = \exp(-i\theta) \cdot \bar{\mathbf{z}}\),
hence \(\bar{\mathbf{z}}\) is the eigenvector of \(\Phi\)
corresponding to the eigenvalue \(\exp(-i\theta)\). The vectors \(\mathbf{z}\)
and \(\bar{\mathbf{z}}\) are linearly independent hence the vectors
\(\mathbf{u} = \dfrac{\mathbf{z} + \bar{\mathbf{z}}}{2}\) and
\(\mathbf{v} = \dfrac{\mathbf{z} - \bar{\mathbf{z}}}{2i}\)
are linearly independent too and can be
taken as the base of the space \(\mathbf{R}^{2}\). It follows from
\eqref{eq19:art2} and \eqref{eq20:art2} that the matrix for the operator
\(\varphi\) in the base \(\mathbf{u}, \mathbf{v}\) is
\begin{equation}
 \label{eq21:art2}
A^{\prime\prime} =
\begin{pmatrix}
\cos\theta & -\sin\theta \\
\sin\theta & \cos\theta
\end{pmatrix}.
\end{equation}

This matrix describes a rotation of the plane on an angle \(\theta \) if and only
if the basis \(\mathbf{u}, \mathbf{v}\) is orthonormal, but for an operator
with the complex spectrum it is not true in general.

Let us look at two co-ordinate systems on the plane. One of them has the
\(x\)-axis and the \(y\)-axis corresponding
the base \(\mathbf{i}, \mathbf{j}\) and another co-ordinate system has the
\(x^{\prime}\)-axis and the \(y^{\prime}\)-axis corresponding the base \(\mathbf{u}, \mathbf{v}\). Let
\(x^{\prime}\) and \(y^{\prime}\) are coordinates of vector \(\mathbf{x}\) in the base \(\mathbf{u}, \mathbf{v}\),
then according to \eqref{eq21:art2}, coordinates of vector
\(\varphi(\mathbf{x})\) in this base, \(\cos\theta \cdot x^{\prime} - \sin\theta \cdot y^{\prime}\)
and \(\sin\theta \cdot x^{\prime} + \cos\theta \cdot y^{\prime}\), satisfy the condition
\((\cos\theta \cdot x^{\prime} - \sin\theta \cdot y^{\prime})^{2} + (\sin\theta \cdot x^{\prime} + \cos\theta \cdot y^{\prime})^{2} =
(x^{\prime})^{2} + (y^{\prime})^{2}\).

It means, that under the action of \(\varphi\), each point \(M_0(x_0^{\prime},y_0^{\prime})\)
remains on a curve \(\Gamma \)
\begin{equation}
 \label{eq22:art2}
\Gamma\colon(x^{\prime})^{2} + (y^{\prime})^{2} = (x_0^{\prime})^{2} + (y_0^{\prime})^{2}.
\end{equation}

Let \(P=
\begin{pmatrix}
p & q \\
s & t
\end{pmatrix}\)
be the change-of-basis matrix (from the base \(\mathbf{u},\mathbf{v}\) to
the base \(\mathbf{i},\mathbf{j}\)), then
\begin{equation}
 \label{eq23:art2}
\begin{pmatrix}
x \\
y
\end{pmatrix} = P^{-1}
\begin{pmatrix}
x^{\prime} \\
y^{\prime}
\end{pmatrix}
\quad \text{and} \quad
\begin{pmatrix}
x^{\prime} \\
y^{\prime}
\end{pmatrix} = P
\begin{pmatrix}
x \\
y
\end{pmatrix}.
\end{equation}

Now we are able to get the equation of curve \(\Gamma\) in the
Cartesian coordinate system. \(\Gamma\colon (x^{\prime})^2+(y^{\prime})^2 =
(px+qy)^2+(sx+ty)^2 = (p^2+s^2)x^2 + 2(pq+st)xy + (q^2+t^2)y^2=r^2\) where
\(r^2 = (x_0^{\prime})^2+(y_0^{\prime})^2\).
The curve \(\Gamma\) is the second-order curve and its quadric quantic is
\(f(x,y) = (p^2+s^2)x^2 + 2(pq+st)xy + (q^2+t^2)y^2\).
\(\mathbf{A}_f=
\begin{pmatrix}
p^2+s^2 & pq+st \\
pq+st & q^2+t^2
\end{pmatrix}\)
is the matrix of this quadric quantic.

The three fundamental invariants \(\Delta\), \(\delta\) and \(S\) determine a second-order
curve up to a motion of the Euclidean plane:
\(\Delta =
\begin{vmatrix}
p^2+s^2 & pq+st & 0 \\
pq+st & q^2+t^2 & 0 \\
 0 & 0 & -r^2
\end{vmatrix}\),
\(\delta =
\begin{vmatrix}
p^{2}+s^{2} & pq+st \\
pq+st & q^{2}+t^{2}
\end{vmatrix}\),
\(S = p^{2}+s^{2} + q^{2}+t^{2}\).

As \(S>0\), \(\delta = (p^{2}+s^{2})(q^{2}+t^{2}) - (pq+st)^{2} = (sq-pt)^{2}
= \det^2 (P)>0\) and \(\Delta = -\delta r^{2}<0\), then \(\Gamma \) is
ellipse. The directions of major axis and minor axes coincide with the
directions of the eigenvectors corresponding to the eigenvalues,
\(\mu^{\prime} > \lambda^{\prime} >0\)
of the matrix \(A_{f}\). So the canonical equation of an ellipse is
\(\dfrac{x^2}{a^2}+\dfrac{y^2}{b^2}=1\), where \(a =
\dfrac{r}{\sqrt{\lambda^{\prime}}}\), \(b = \dfrac{r}{\sqrt{\mu^{\prime}}}\) and irrespective
of \(r\) all these ellipses are similar. Thus under the action of \(\varphi\)
each point of the plane remains on the ellipse.
The trajectory of the points is either finite number of points of ellipse
if the angle \(\theta\) is commensurable with \(\pi\) or the infinite
dense set of points in the opposite case.

\renewcommand{\refname}{Literature:}

\vspace{1cm}
\par\medskip\noindent
email: ibusjatskaja@hse.ru, yukochetkov@hse.ru

\end{document}